\documentclass[12pt]{article}

\usepackage[latin1]{inputenc}                      
\usepackage[T1]{fontenc}                           
\usepackage[english]{babel}                        

\usepackage[a4paper,dvips,height=240mm,width=165mm]{geometry}

\usepackage{amsmath,amsfonts,bbm,stmaryrd}
\usepackage{graphicx}

\usepackage{theorem}{\theorembodyfont{\slshape}
\newtheorem{theorem}{Theorem}[section]             
\newtheorem{proposition}[theorem]{Proposition}    
\newtheorem{remark}[theorem]{Remark}              
\newtheorem{lemma}[theorem]{Lemma}                 
}

\makeatletter                                                                             
  
\@addtoreset{equation}{section} 						
\makeatother  										

\def\Proof{{\smallskip\noindent{\em Proof. }}}     
\def\endProof{{\hfill$\Box$\medskip\noindent}}    

\newcommand\R{{\mathbb R}}                         
\newcommand\Y{{\mathbb Y}}                         
\renewcommand\H{{\mathbb H}}                       
\renewcommand\P{{\mathbb P}}                       
\newcommand\B{{\mathbb B}}                         
\newcommand\U{{\mathbb U}}                         
\newcommand\V{{\mathbb V}}                         
\let\div=\relax\DeclareMathOperator{\div}{div}     
\newcommand{\indicatrice}{\mathbbm{1}}             
\newcommand{\semicolon}{\,;\,}                     
\let\theta=\vartheta                               
\let\tilde=\widetilde                              
\newcommand{\definition}{\underset{\text{def}}{=}} 

\title{On the localization of the magnetic and the
velocity\\ fields in the equations of magnetohydrodynamics}

\author{Lorenzo Brandolese$^\ast$ and 
	Fran\c cois Vigneron$^\dagger$}

\newcommand{\makeadress}{
$^\ast$
\begin{minipage}[t]{.4\textwidth}
\begin{flushleft}
Lorenzo Brandolese\\
Institut Camille Jordan.\\
Universit\'e Claude Bernard Lyon 1.\\
21 avenue Claude Bernard\\
F-69622 Villeurbanne Cedex\\
FRANCE\\
{\tt brandolese@math.univ-lyon1.fr}
\end{flushleft}
\end{minipage}
\hfill
$^\dagger$
\begin{minipage}[t]{.45\textwidth}
\begin{flushleft}
Fran\c cois Vigneron\\
Centre de Math\'ematiques L.~Schwartz.\\
U.M.R. 7640 du C.N.R.S.\\
Ecole polytechnique.\\
F-91128 Palaiseau Cedex\\
FRANCE\\
{\tt francois.vigneron@normalesup.org}
\end{flushleft}
\end{minipage}
}

\begin{document}

\maketitle

\begin{abstract}
We study the behavior at infinity, with respect to the space variable,
of solutions to the magnetohydrodynamics equations in $\R^d$.
We prove that if the initial magnetic field decays sufficiently fast,
then the plasma flow behaves as a solution of the free nonstationnary
Navier--Stokes equations when $|x|\to +\infty$, and that the magnetic field
will govern the decay of the plasma, if it is poorly localized at the
beginning of the evolution. Our main tools are boundedness criteria for
convolution operators in weighted spaces.
\end{abstract}

\section{Introduction}
The magnetohydrodynamics equations are a well-known model in plasma physics,
describing the interactions between a magnetic field and a fluid made of
moving electrically charged particles.
A common example of an application of this model is the design of tokamaks:
the purpose of these machines is to confine a plasma in a region, with a
density and a temperature large enough to entertain thermonuclear fusion
reactions. This can be achieved, at least during a small time interval,
by applying strong magnetic fields.
We refer to \cite{Pri82} for other applications of this model, in particular
to the study of the dynamics of the solar corona.

In non-dimensional form, the magnetohydrodynamics equations can be written
in the following way:
\begin{equation*}
\left\{
\begin{aligned}
&\frac{\partial u}{\partial t} + (u\cdot\nabla)u-S(B\cdot\nabla)B
+\nabla\left(p+\frac{S}{2}|B|^2\right)=\frac{1}{R_e}\Delta u \\
&\frac{\partial B}{\partial t} +(u\cdot\nabla)B-(B\cdot\nabla)u=
\frac{1}{R_m}\Delta B \\
&\div u=\div B=0\\
&u(0)=u_0\quad\text{and}\quad B(0)=B_0.
\end{aligned}
\right.
\eqno{\mbox{(MHD)}}
\end{equation*}
Here the unknowns are the velocity field~$u$ of the fluid, the pressure~$p$
and the magnetic field~$B$, all defined in $\R^d$ $(d\ge2)$.
The positive constants $R_e$ and $R_m$ are respectively the Reynolds number
and the magnetic Reynolds number; moreover $S=M^2/(R_eR_m)$, where $M$ is 
the Hartman number.
After rescaling $u$ and $B$, we can assume that $S=R_e=1$. 
With minor loss of generality, from now on we shall also assume that $R_m=1$.
All the results however remain valid in the general case with
simple modifications in the constants.

In the particular case $B\equiv 0$, the system (MHD) reduces to the celebrated
Navier--Stokes equations. Just as in this particular case, global weak solutions
to (MHD) do exist, but their unicity, as well their smoothness in the case of
smooth data, remains an open problem for $d\ge3$.
Partial regularity results, which provide bounds of the Hausdorff dimension of
the possible singular set of weak solutions, have been obtained in \cite{HeX05i}.
Constantin and Fefferman's theory \cite{CosF93} relating the regularity of
the flow to the directions of the vorticity has been extended to
magnetohydrodynamics in~\cite{HeX05ii}.
A construction of forward selfsimilar solutions is given in \cite{HeX05iii}, where the nonexistence
of backward selfsimilar solutions is also discussed.
Moreover, the asymptotic behavior of the solutions for $t\to +\infty$ is quite
well understood: for example, \cite{SchSS} provides the optimal decay rates
of the $L^2$ norm of $u$ and $B$ for a large class of flows.

On the other hand, nothing seems to have been done to study the decay of
solutions of (MHD) with respect to the \emph{space variable}. In this paper,
motivated by recent results obtained by several authors for the Navier--Stokes
equations (see, {\it e.g.\/}, \cite{Bra04iii}, \cite{BraM02}, \cite{HeM03},
\cite{Mi00} and \cite{Vig05}), we would like to describe in which way the
presence of the magnetic field affects the spatial localization of the
velocity field.

\paragraph{Definitions and notations.}

We start by introducing the notion of decay rate at infinity in a {\emph weak sense\/},
which generalizes the usual notion of pointwise decay rate in the framework
of locally square integrable functions.
A simple motivation is that the $L^2_{\text{loc}}$ regularity is the minimal
one for which the system (MHD) makes sense.

\begin{enumerate}
\item
Let  $f\in L^2_{\text{loc}}(\R^d)$.
We define the  \emph{$L^2$ decay rate as $|x|\to +\infty$} of~$f$, as
\begin{equation}
\label{distr-decay}
\eta(f)=\sup\biggl\{\eta\in\R\;;\;  
\lim_{R\to +\infty} R^{2\eta} \int_{1\le|x|\le 2} |f(Rx)|^2  \, dx =0
\biggr\}.
\end{equation}
If $\eta=\eta(f)$ is finite then we will write 
$f\stackrel{L^2}{\sim} |x|^{-\eta}$ when $|x|\to +\infty$.
On the other hand, when we write $f\stackrel{L^2}{=}{\cal
O}(|x|^{-\eta})$ when $|x|\to +\infty$,
we mean that $\eta(f)\ge \eta$.
Of course, any measurable function such that $|f(x)|\le C(1+|x|)^{-\eta}$
satisfies $f\stackrel{L^2}{=}{\cal O}(|x|^{-\eta})$ when $|x|\to +\infty$.

\item
For $a\in[1,+\infty]$ and $\alpha\in\R$,
the space $L^a_\alpha(\R^d)$ is the Banach space normed by
\begin{subequations}
\begin{equation}
\|f\|_{L^a_\alpha}=
\biggl(\int_{\R^d} |f(x)|^a(1+|x|)^{a\alpha}\,dx\biggr)^{1/a}
\qquad \hbox{if $1\le a<+\infty$}
\end{equation}
and, if  $a=+\infty$, by
\begin{equation}
\|f\|_{L^\infty_\alpha}= \underset{x\in\R^d}{\operatorname{ess\enspace sup}}
\enspace|f(x)|(1+|x|)^\alpha.
\end{equation}
\end{subequations}

From the \emph{localization point of view}
the two spaces $L^a_\alpha(\R^d)$ and $L^b_\beta(\R^d)$ must be considered
as equivalent, when 
$$\alpha+\frac{d}{a}=\beta+\frac{d}{b}.$$
Indeed, if $ f\in L^a_\alpha(\R^d)$ and $a\geq2$, then 
$f\stackrel{L^2}{=}{\cal O}\bigl(|x|^{-(\alpha+d/a)}\bigr)$
when $|x|\to +\infty$.
H\"older inequality implies that
\begin{equation}\label{inclusion}
L^a_\alpha \subset L^b_\beta
\end{equation}
whenever $\alpha+d/a > \beta+ d/b$ and $a\geq b$. It also implies that
\begin{equation}\label{def2_eta}
\eta(f) = \sup\left\{ \alpha+\frac{d}{a} \,;\, a\geq 2 \enspace\text{and}
\enspace f\in L^a_\alpha \right\}
\end{equation}
for any $f\in L^2_{\text{loc}}(\R^d)$.

\end{enumerate}

\noindent
We shall use the following additional notations :

\begin{enumerate}
\setcounter{enumi}{2}

\item
If  $A$ and $B$ are two expressions containing a parameter $\alpha$,
then when we write
$$ A\le B-\varepsilon_\alpha,$$
we mean that
$A\le B$ if $\alpha=0$ and $A<B$ if $\alpha\not=0$.
We shall also often write expressions of the form $A\le B-\varepsilon_{1/a}$
meaning that the inequality must be strict for finite $a$
and can be large when $a=+\infty$.

\item
The positive part of a real number will be denoted
by $(\cdot )^+=\max\{\cdot,0\}$.

\end{enumerate}

\paragraph{Main results.}

We are concerned  with the persistence problem of the spatial localization
of the magnetic and the velocity fields.
Our main results (Theorem~\ref{main-thm} and \ref{theorem1} below) aim to
answer the following questions. Consider  a localization condition like 
\begin{equation}  
\label{assumption-data}
(u_0,B_0)\in L^{p_0}_{\theta_0}(\R^d)\times L^{p_1}_{\theta_1}(\R^d).
\end{equation}
Will the unique solution of (MHD) preserve such a condition in some future
time interval ? Depending on the parameters, the answer can be positive
or negative. In case of a negative answer, can we still ensure that the
spatial localization of the solution is conserved {\it in the weak sense\/} ?
In other words, we would like to know whether
$$u(t)\stackrel{L^2}{=}{\cal O}\bigl(|x|^{-(\theta_0+d/p_0)}\bigr)
\quad\hbox{and}\quad
 B(t)\stackrel{L^2}{=}{\cal O}\bigl(|x|^{-(\theta_1+d/p_1)}\bigr)\qquad
 \hbox{when $|x|\to +\infty$}.$$
Again, this condition may be conserved, or instantaneously break down.

\medskip 
We will prove the following:

\begin{theorem}
\label{main-thm}
Let $u_0\in L^{p_0}_{\theta_0}(\R^d)$, $B_0\in L^{p_1}_{\theta_1}(\R^d)$ be
two divergence-free vector fields in $\R^d$ ($d\ge2$).
Assume that
\begin{subequations}\label{hyp-thm}
\begin{equation}\label{indices-thm}
\left\{\begin{gathered}
\theta_0\geq 0\\ d<p_0\leq +\infty
\end{gathered}\right.
\qquad\text{and}\qquad
\left\{\begin{gathered}
\theta_1\geq0 \\ d<p_1\leq +\infty.
\end{gathered}\right.
\end{equation}
Let us also assume that
\begin{equation}\label{eta-thm}
\delta + \varepsilon_\delta \leq \eta_0\le \min\big\{ d+1\semicolon 2\eta_1-\delta \big\},
\end{equation}
\end{subequations}
with $\eta_0=\theta_0+d/p_0$, $\eta_1=\theta_1+d/p_1$ and
$\delta=\left(\frac{2d}{p_1}-1\right)^+$\!.
Finally, define
$p_0^\ast=\min\{p_0\semicolon\frac{d}{\delta}-\varepsilon_{\delta}\}$.
\par\smallskip
Then there exists $T>0$ and a unique mild solution $(u,B)$ of (MHD) in
${\cal C}([0,T];L^{p_0^\ast}\times L^{p_1})$.
This solution satisfies
\begin{equation}\label{stability-thm}
u(t)\overset{L^2}{=}{\cal O}\bigl(|x|^{-\eta_0}\bigr)
\quad\text{and}\quad
 B(t)\overset{L^2}{=}{\cal O}\bigl(|x|^{-\eta_1}\bigr)
 \qquad\text{when $|x|\to +\infty$}.
\end{equation}
If $d=2$, the time $T$ can be arbitrarily large.
\par\smallskip\pagebreak[1]
Moreover, if $(u_0,B_0)$ also belongs to
$L^{\tilde{p_0}}_{\tilde\theta_0}\times L^{\tilde{p_1}}_{\tilde\theta_1}$,
with the corresponding indices satisfying assumptions~\eqref{hyp-thm},
then the lifetimes in $L^{p_0^\ast}\times L^{p_1}$ and
$L^{\tilde{p_0}^{\!\ast}}\times L^{\tilde{p_1}}$ agree and both maximal solutions
are actually the same one.
\end{theorem}

\bigskip
Next we discuss the optimality of the above restrictions.
Such restrictions are of two kinds: 
there are a few conditions related to the well-posedness of
the system, and a condition (namely, the upper bound for $\eta_0$ in  \eqref{eta-thm})
which is related to the spatial localization of the solution.
Here, we will only focus on this condition.
The following theorem implies that the restriction $\eta_0\leq d+1$ is sharp.
We expect that the other restriction is also sharp, or at least
that $\eta_0\leq 2\eta_1$ for stable weak solutions. But we were not
able to prove such a result.

\begin{theorem}\label{thm_optimal}
Let $(u,B)\in {\cal C}([0,T];L^2(\R^d)\times L^2(\R^d))$ a solution to (MHD)
such that
\begin{subequations}
\begin{align}\label{optimality_u}
\sup_{t\in[0,T]} |u(t,x)| & \overset{L^2}{=} {\cal
O}\bigl(|x|^{-(d+1+\varepsilon)}\bigr)\\
\label{optimality_B} \text{and}\qquad
\sup_{t\in[0,T]} |B(t,x)| & \overset{L^2}{=} {\cal
O}\bigl(|x|^{-(d+1+\varepsilon)/2}\bigr)
\end{align}
\end{subequations}
for some $\varepsilon>0$. Then, for all $t\in[0,T]$, there exists a  constant
$C(t)\geq0$ such that the components of $u(t)$ and $B(t)$ satisfy the
following integral identity :
\begin{equation}
\label{orthogonality}
\int_{\R^d} (u^ju^k-B^jB^k)(t,x)\,dx= \delta_{j,k}\,C(t), \qquad (j,k=1,\ldots,d)
\end{equation}
with $\delta_{j,k}=1$ if $j=k$ and~$\delta_{j,k}=0$ otherwise.
\end{theorem}
By Theorem~\ref{theorem1} below, 
condition \eqref{optimality_B} will be fulfilled as soon as $u_0$
and $B_0$ belong to $L^p_\theta(\R^d)$,
with $p>d$ and $\theta+\frac{d}{p}=(d+1+\varepsilon)/2$,
for some $\epsilon>0$.
This means that if we start with a well localized initial datum $(u_0,B_0)$,
but such that \eqref{orthogonality} does not hold for $t=0$,
then condition \eqref{optimality_u} must brake down.

\bigskip
On the other hand, the integral identities \eqref{orthogonality} are obviously unstable.
Neverthless, in section~\ref{section-spreading}
we shall see that a class of exceptional solutions satisfying \eqref{orthogonality} does exist.
Inside this class, one can exhibit solutions such that $u$ decays much faster
than in the generic case.

\begin{figure}[!htb]\label{figure}\noindent\begin{tabular}{cc}
 \includegraphics{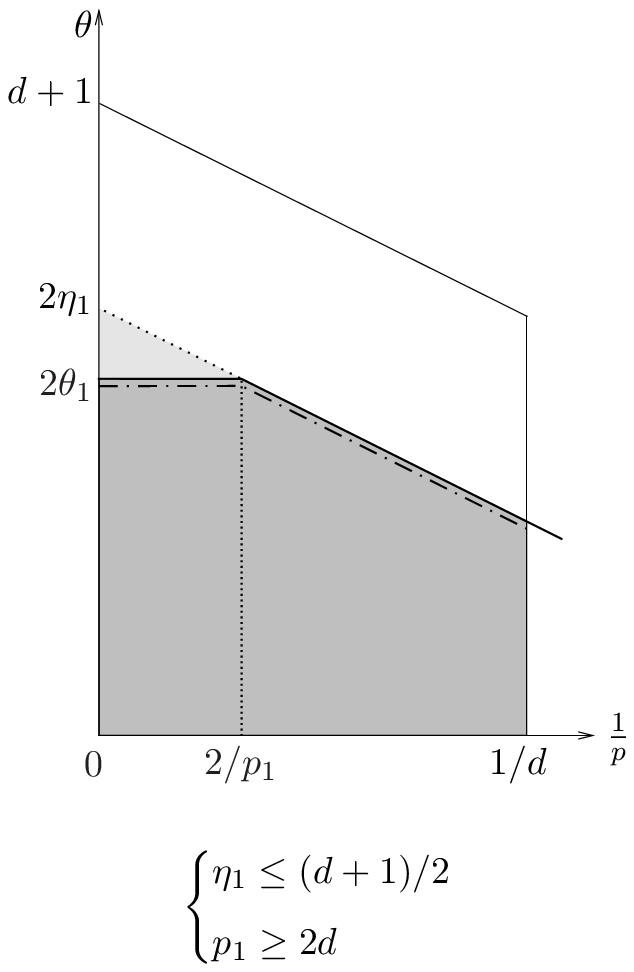} &
 \includegraphics{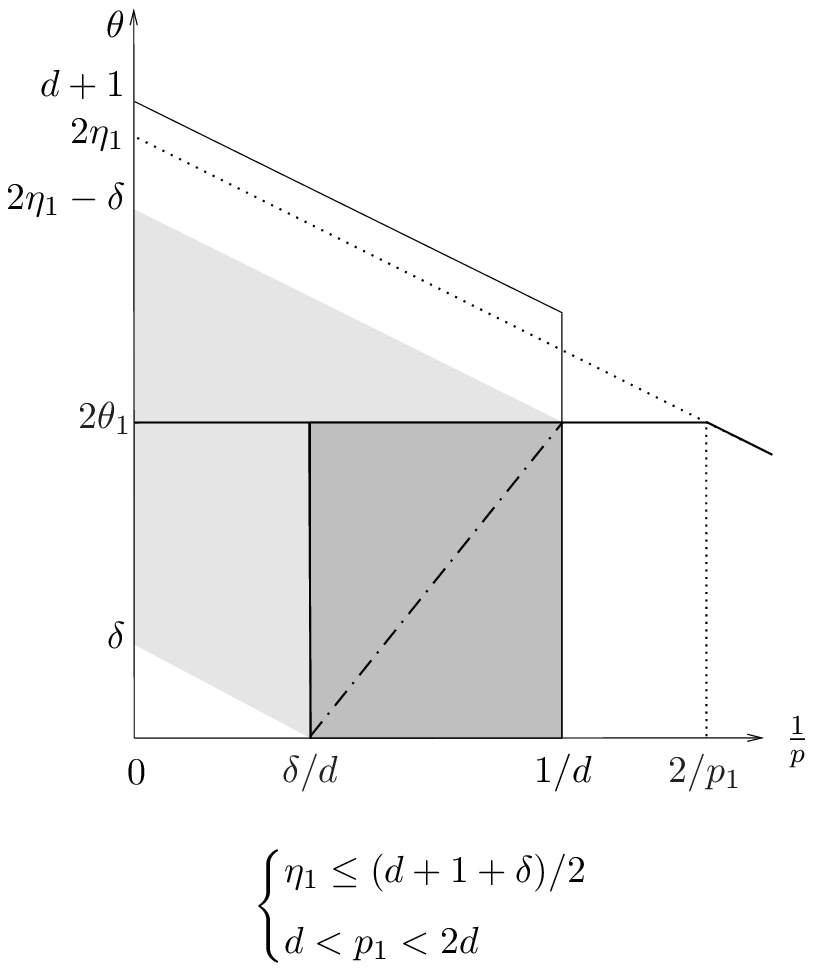} \\[1.5em]
 \includegraphics{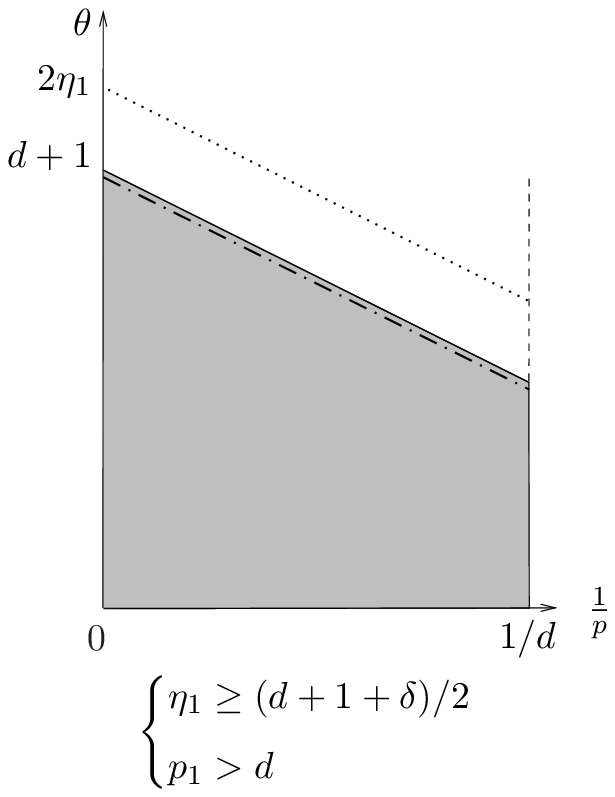}  &
\begin{minipage}[b]{.45\textwidth}
{\bf Fig.1} \quad
The figures show the admissible values for $(p_0,\theta_0)$
allowing \eqref{stability-thm} to hold, once $(p_1,\theta_1)$ is given
({\sl all gray regions}).
\par
\medskip
{\sc Above} : Slowly decaying magnetic field. The results depends slightly on
the regularity of $B$ through $\delta=\left(\frac{2d}{p_1}-1\right)^+$\!.
\par
\medskip
{\sc Down-Left} : Fast decaying magnetic field. The velocity field
behaves at infinity as the solution of Navier--Stokes equations with the
same initial datum $u_0$ (see \cite{Vig05}).
\par
\bigskip
The {\sl dark gray} regions correspond to initial data for wich
we will prove in addition that $u\in L^\infty([0,T];L^{p_0}_{\theta_0})$.
The {\sl dash-dotted lines} illustrate the barriers used in the
proof of~\S\ref{proof-main-thm}.\vspace*{1.5ex}
\end{minipage}
\end{tabular}
\end{figure}

\paragraph{Physical interpretation of Theorem~\ref{main-thm}.}

This theorem reinforces mathematically some facts
that can be observed in the applications.
Three conclusions can be drawn:
\begin{enumerate}
\item Any spatial localization assumption on the magnetic field will be
conserved by the flow. Indeed, the $L^2$~decay rate  $\eta_1$
can be arbitrarily large.
The spatial localization of the velocity field is also conserved,
but there are some limitations to this property.
\item
For poorly localized magnetic fields (namely $\eta_1\leq(d+1+\delta)/2$),
the behavior of $u$ when $|x|\to +\infty$ is governed by the decay of the magnetic
field. As $0\leq\delta<1$ in \eqref{eta-thm}, the maximal $L^2$~decay rate of $u$ that can be
conserved by the flow exceeds $2\eta_1-1$.
When $p_1\ge 2d$, one has~$\delta=0$ and this rate is improved up to
twice that of~$B_0$. The pathological lower bound on
$\eta_0$ disappears too.
Roughly speaking, requiring $p_1$ to be larger (for a given $L^2$~decay rate $\eta_1=\theta_1+d/p_1$
of the magnetic field) means that the behavior at infinity of $B_0$ is closer
and closer to that of a function that decays as~$|x|^{-\eta_1}$, in the usual
pointwise sense.
\item For sufficiently fast decaying magnetic fields,
the decay of  $u$ is not affected by $B$, but is provided by the fundamental
laws of hydrodynamics.
The reason is the following: 
for magnetic fields such that $\eta_1\ge (d+1+\delta)/2$, our limitations
on the $L^2$ decay rate at infinity of the velocity field  \eqref{eta-thm}
boil down to the only restriction $\eta_0\leq d+1$.
This is exactly the same restriction that appears for the Navier--Stokes 
equations. Indeed, we know from F.~Vigneron's result \cite{Vig05}
that the mild solution of the Navier--Stokes equations remains in 
$L^{p_0}_{\theta_0}(\R^d)$ if the initial velocity belongs to such space
and 
\begin{equation*}
\theta_0+d/p_0\le d+1-\varepsilon_{1/p_0}.
\end{equation*}
This condition in known to be sharp. One may notice however that,
thanks to \eqref{def2_eta}, the equality case is possible
even if~$p_0<+\infty$, provided that stability is asserted as
in \eqref{stability-thm}.
\end{enumerate}
A more physical explanation for the above conclusions is the following\footnote
{This explanation was suggested to us by the Referee.}.
The induction equation means that the magnetic field lines
are transported by the flow while simultaneously undergoing resistive
diffusion.
This transport-diffusion process guarantees that, where the velocity vanishes,
the magnetic field
will not spatially spread out during small time intervals,
since the mechanism of diffusion is quite slow.
As for the fluid flow, the magnetic field acts upon it only through the Lorentz force:
whenever this disappears the velocity acts in a purely Navier--Stokes way;
thus, the spatial spreading of the initial velocity is essentially governed by 
the competition between diffusion, whose effect is important only for large time,
and incompressibility, that immediately prevents the flow from remaining too localized.

\paragraph{Stability in weighted spaces.}

Conclusion \eqref{stability-thm} does not mean that
$$(u,B)\in L^\infty\left([0,T] \semicolon L^{p_0}_{\theta_0}\times
    L^{p_1}_{\theta_1}\right).$$
Actually, we do not know if this property holds
when $u_0\in L^{p_0}_{\theta_0}$ and $(p_0,\theta_0)$
is in the light-gray regions of~Fig.1.
However, if $(p_0,\theta_0)$ is in a dark-gray region, then such property
does hold.
This is essentially the statement of our next theorem.
It extends to the case of non-vanishing magnetic fields, the result
established in \cite{Vig05} for the Navier--Stokes equations.

\begin{theorem}\label{theorem1}
Let $u_0\in L^{p_0}_{\theta_0}(\R^d)$, $B_0\in L^{p_1}_{\theta_1}(\R^d)$ be
two divergence-free vector fields in $\R^d$ ($d\ge2$).
\begin{subequations}
Assume that $\theta_0,\theta_1\geq 0$, $d<p_0\le +\infty$ and
\begin{equation}\label{hypo}
\frac{2}{p_1}<\frac{1}{p_0}+\frac{1}{d}.
\end{equation}
Then there exist $T>0$ (if $d=2$, one may take~$T=+\infty$)
and a unique mild solution of (MHD)
\begin{equation}\label{mild-solution_Lp}
(u,B)  \in \mathcal{C}\left([0,T];L^{p_0}\times L^{p_1}\right).
\end{equation}
\end{subequations}
\begin{subequations}
If, in addition, the decay rates of $u_0$ and $B_0$ defined by
$\eta_0=\theta_0+d/p_0$ and
$\eta_1=\theta_1+d/p_1$ satisfy
\begin{equation}
\label{local}
\eta_0\le\min
\biggl\{ d+1-\varepsilon_{1/p_0}\semicolon 2\eta_1-\varepsilon_{2\theta_1-\theta_0}\semicolon
             2\eta_1+\frac{d}{p_0}-\frac{2d}{p_1} \biggr\},
\end{equation}
then we have more precisely
\begin{equation}\label{mild-solution}
(u,B)\in {\cal C}\bigl([0,T];L^{p_0}_{\theta_0}\times L^{p_1}_{\theta_1}\bigr).
\end{equation}
Moreover, if $(u_0,B_0)$ also belongs to
$L^{\tilde{p_0}}_{\tilde\theta_0}\times L^{\tilde{p_1}}_{\tilde\theta_1}$,
with new indices again satisfying~\eqref{hypo} and~\eqref{local},
then the lifetimes in $L^{p_0}_{\theta_0}\times L^{p_1}_{\theta_1}$ and
$L^{\tilde{p_0}}_{\tilde \theta_0}\times L^{\tilde{p_1}}_{\tilde \theta_1}$
are the same and both maximal solutions agree.
\end{subequations}
\end{theorem}

The assumption \eqref{hypo} is not really related to spatial localization
problems, but rather to well-posedness issues of the equations, and in
particular, to the invariance of the equation under the natural scaling
$$u_\lambda(t,x)=\lambda u(\lambda^2t,\lambda x),\qquad
B_\lambda(t,x)=\lambda B(\lambda^2t,\lambda x) \qquad (\lambda>0).$$
We expect that Theorem~\ref{theorem1} remains true in limit cases~$p=d$,
or $\frac{2d}{p_1}=\frac{1}{p_1}+\frac{1}{d}$
(with several modifications in the proof).
We did not treat these limit cases since they would require 
Kato's two-norm approach 
for proving the boundedness of the operators involved, as described in
\cite[chap.~3]{Can95} or \cite{CanK04} for the Navier--Stokes equations.
The proof would be more complicated,
without providing any substantial clarification of the spatial localization problem.

Let us also observe that one could replace the weights $(1+|x|)^\theta$
with homogeneous weights. But in this case the conditions to be imposed 
on the parameters would be much more restrictive, {\it e.g.}
$$\theta+\frac{d}{p}<1.$$
Again, this would not help to understand the spatial localization of the
fields.

\paragraph{Main methods and organization of the paper.}
We shall first prove Theorem~\ref{theorem1}
and later deduce Theorem~\ref{main-thm} as a corollary
of the natural embedding~\eqref{inclusion} between weighted spaces.
The idea consists in observing that
the assumptions \eqref{hyp-thm}, together with the
inclusion \eqref{inclusion}, ensure that the initial datum belongs to the product of
two larger Lebesgue spaces, in which we can prove the existence and
uniqueness of a mild solution.

Our proof of Theorem~\ref{theorem1} consists in applying
the contraction mapping principle to the integral form of (MHD),
in a suitable ball of the
space~${\cal C}([0,T],L^{p_0}_{\theta_0}\times L^{p_1}_{\theta_1})$.
This is why we refer to $(u,B)$ as a mild solution.
The only difficulty is establishing
the bicontinuity of the bilinear operator involved.

For small values of $\eta_0$,
the bicontinuity would be a straightforward consequence of the well-known Young convolution
inequality in weighted Lebesgue 
spaces (recalled in \cite[\S2.2]{Vig05}).
But this argument does not go through 
when $\eta_0$ is close to the upper bound of~\eqref{local},
since the kernel of the operator governing the evolution
of the velocity field decays too slowly at infinity.
In this case, the proof requires more careful estimates.
The main one is given by Proposition~\ref{proposition1} below.

Several generalizations of the weighted convolution inequalities are known
(see, {\it e.g.\/}, the recent boundedness criterion for asymmetric kernel operators
\cite[\S2.3]{Vig05}, which applies to Navier--Stokes).
However, we could not deduce the bicontinuity of the bilinear operator
by applying directly any known inequality,
unless we put additional artificial restrictions on the parameters.

The main issue with the spatial localization of magnetohydrodynamics fields
is that the system cannot be treated as a scalar equation. When dealing with
the Navier--Stokes system, one may often reduce the problem to a single
equation, because all the components of the kernels of the Navier--Stokes
operators satisfy the same estimates. This is no longer true for
(MHD). In the following, we shall derive sharp bounds for
the magnetohydrodynamics kernels and take advantage of the fact that a few
components decay much faster than the others.

\bigskip
This paper is organized as follows.
Section~\ref{section-generalities}
contains  some generalities on magnetohydrodynamics.
In Section~\ref{section-convolution} we study the boundedness of convolution
operators in weighted spaces.
We use these results in Section~\ref{section-lifetime}, proving first
the local existence of a unique solution in weighted
spaces \eqref{mild-solution}, then the fact that lifetimes do not depend
on the choice of the indices.
Then we deduce Theorem~\ref{main-thm} as a corollary.

Theorem~\ref{thm_optimal} will be proved in
Section~\ref{section-spreading},
using a Fourier transform method developed in \cite{BraM02}.
Section~\ref{section-spreading} also contains the
description of a method for obtaining special solutions,
such that the velocity field is more localized than in \eqref{eta-thm}.
Those solutions are however unstable.

\begin{remark}
\begin{rm}
When we deal with the space
${\cal C}([0,T];L^{p_0}_{\theta_0}\times L^{p_1}_{\theta_1})$,
with $p_0=+\infty$ or $p_1=+\infty$, the continuity at $t=0$ must be understood in the
weak sense, as is usually done in nonseparable spaces.
\end{rm}
\end{remark}


\section{The integral form of the equations}
\label{section-generalities}

Let $\P$ be the Leray-Hopf projector onto the divergence-free vector field,
defined by $$\P f=f-\nabla\Delta^{-1}(\div f).$$
Applying $\P$  to the first equation of (MHD) and then the Duhamel formula,
we obtain the integral equations
\begin{equation*}
\left\{\begin{aligned}
&u(t)=e^{t\Delta}u_0 - \int_0^t e^{(t-s)\Delta}\:\P\div(u\otimes u - B\otimes B)(s)\,ds\\
&B(t)=e^{t\Delta}B_0 - \int_0^t e^{(t-s)\Delta}\:\div(u\otimes B - B\otimes u)(s)\,ds\\
&\div u_0=\div B_0=0
\end{aligned}
\right.
\eqno\mbox{(IE)}
\end{equation*}
where $e^{t\Delta}$ is the heat semigroup (recall that the Reynolds numbers
and the Hartman numbers have been set equal to $1$).
The semigroup method that we use in this paper to solve (IE)
provides mild solutions of (MHD) that are in fact smooth for strictly
positive~$t$.

We denote respectively by $F_{j,h}^k(t,x)$ and
$G_{j,h}^k(t,x)$ $(j,h,k=1,\ldots,d)$
the components of the kernels of the matricial operators
$e^{t\Delta}\P\nabla$ and $e^{t\Delta}\nabla$.
Thus,
\begin{equation}
\label{symbol}
\widehat{F}_{j,h}^k(\xi,t)=
  e^{-t|\xi|^2}\xi_h(\delta_{j,k}-\xi_j\xi_k|\xi|^{-2}).
\end{equation}
This expression of the symbol allows us to see that
\begin{subequations}
\begin{align}
\label{decay-F}
\begin{split}
\mbox{}&F(t,x)=t^{-(d+1)/2}\:\Phi(x/\sqrt t),\\
\text{with}\quad&|\Phi(x)|\le C(1+|x|)^{-(d+1)}.
\end{split}
\intertext{This low decay rate of $\Phi$ is due to the fact that
$F(t,\cdot)\not\in L^1_1(\R^d)$ ; otherwise $\widehat F(t,\cdot)$ would be a
${\cal C}^1$ function on $\R^d$. On the other hand,}
\label{decay-G}
\begin{split}
\mbox{}&G(t,x)=t^{-(d+1)/2}\:\Psi(x/\sqrt t),\\
\text{with}\quad&\Psi\in{\cal S}(\R^d) \quad\text{(the Schwartz class)}.
\end{split}
\end{align}
\end{subequations}

Let us introduce the bilinear operators on $\R^d$-vector fields $\U$ and
$\B$ whose $k$\textsuperscript{th} component is
\begin{eqnarray*}
\U^k(f,g)(t,x)&=&\sum_{j,h}\int_0^t F_{j,h}^k(t-s)\ast\left(f^j\otimes
g^h\right)(s)\:ds\\
\B^k(f,g)(t,x)&=&\sum_{j,h}\int_0^t G_{j,h}^k(t-s)\ast\left(f^j\otimes
g^h\right)(s)\:ds,
\end{eqnarray*}
and the bilinear operator $\V=(\V_1,\V_2)$ on $\R^{2d}$-vector fields
$v=(v_1,v_2)$ defined by
\begin{eqnarray*}
\V_1(v,w) &=& \U(v_1,w_1)-\U(v_2,w_2)\\
\V_2(v,w) &=& \B(v_1,w_2)-\B(v_2,w_1).
\end{eqnarray*}
Here and below, for $v\in\R^{2d}$, we denote by $v_1$ the first $d$ components
and by $v_2$ the last $d$ components.

\medskip

With these notations and setting $v=(u,B)$, $v_0=(u_0,B_0)$, the system (IE)
can be rewritten as
\begin{equation}
\label{mhd}
v=e^{t\Delta}v_0- \V(v,v).
\end{equation}

As it is well known (we refer, {\it e.g.\/}, to \cite[Lemma 1.2.6]{Can95}), if $X$ is a
Banach space, then for solving an equation like \eqref{mhd} one just needs
to check that
\begin{subequations}
\begin{equation}
\label{linear}
e^{t\Delta}v_0\in {\cal C}([0,T];X)
\end{equation}
and 
\begin{equation}
\label{bicontinuity}
\V\,:\, {\cal C}([0,T];X)\times {\cal C}([0,T];X) \to {\cal C}([0,T];X),
\end{equation}
\end{subequations}
with the operator norm of $\V$ tending to~$0$ as $T\to0$.
Then the existence of a solution $v\in {\cal C}([0,T];X)$ is ensured,
at least for $T>0$ small enough.

\medskip
In order to prove Theorem~\ref{theorem1} we shall take $X=L^{p_0}_{\theta_0}\times L^{p_1}_{\theta_1}$.
In this setting, condition \eqref{linear}, the unicity and the continuity of
the solution with respect to the time variable are all straightforward.
Therefore, our attention will now be exclusively devoted to the more subtle
problem of the bicontinuity of~$\V$
in~$L^\infty([0,T]; L^{p_0}_{\theta_0}\times L^{p_1}_{\theta_1} )$.

We need three estimates, namely
\begin{subequations}
\begin{align}
\label{Uuu}
\| \U(u,u)(t) \|_{L^{p_0,\vartheta_0}} &\le C_T \| u\|^2_{{\cal C}([0,T],L^{p_0,\vartheta_0})}\\
\label{UBB}
\| \U(B,B)(t) \|_{L^{p_0,\vartheta_0}} &\le C_T \| B\|^2_{{\cal C}([0,T],L^{p_1,\vartheta_1})}\\
\label{BuB}
\| \B(u,B)(t) \|_{L^{p_1,\vartheta_1}} &\le C_T \|u\|_{{\cal C}([0,T],L^{p_0,\vartheta_0})}\| B\|_{{\cal C}([0,T],L^{p_1,\vartheta_1})}
\end{align}
\end{subequations}
for all $0\le t\le T$ and some constant $C_T$ such that $C_T\to0$ as $T\to0$.
These bounds will not rely on the specific structure of the operators $\U$ and $\B$,
but only on the decay properties of their respective kernels:
\begin{equation}
\label{pointwise}
\begin{gathered}
|F(t,x)|\le C (\sqrt t+|x|)^{-(d+1)}\\
|G(t,x)|\le C_N\sqrt t^{N-d-1}(\sqrt t+|x|)^{-N}
\end{gathered}
\end{equation}
for all $N\ge0$.

We start by observing that by H\"older inequality,
\begin{align*}
\|u\otimes u\|_{L^{p_0/2}_{2\theta_0}} 	
	&\le \|u\|_{L^{p_0}_{\theta_0}}^2\\
\|B\otimes B\|_{L^{p_1/2}_{2\theta_1}}		
	&\le \|B\|_{L^{p_1}_{\theta_1}}^2\\
\|u\otimes B\|_{L^{\H(p_0,p_1)}_{\theta_0+\theta_1}}
	&\le \|u\|_{L^{p_0}_{\theta_0}}
 \|B\|_{L^{p_1}_{\theta_1}}
\end{align*}
where $\frac{1}{\H(p_0,p_1)}=\frac{1}{p_0}+\frac{1}{p_1}$
denotes the H\"older exponent (the assumptions of Theorem~\ref{theorem1}
imply that~$p_0,p_1\ge2$).
Set $\lambda=\sqrt t$ and
\begin{equation}
\label{g-lambda,N}
\Gamma_{\lambda}^{N}(x)= (\lambda+|x|)^{-N}.
\end{equation}
Then the only thing that we have to do to obtain \eqref{Uuu}-\eqref{BuB}
is to establish that for all $0<\lambda\le1$:
\begin{subequations}
\begin{align}
\label{uu}
\|\Gamma_{\lambda}^{d+1}*f\|_{L^{p_0}_{\theta_0}} 
	&\le C\lambda^{\sigma_0}
	      \|f\|_{L^{p_0/2}_{2\theta_0}}\,,\\
\label{BB}
\|\Gamma_{\lambda}^{d+1}*f\|_{L^{p_0}_{\theta_0}} 
	&\le C\lambda^{\sigma_0'}
	      \|f\|_{L^{p_1/2}_{2\theta_1}}\,\\
\label{uB}\text{and}\quad
\|\Gamma_{\lambda}^{N}*f\|_{L^{p_1}_{\theta_1}} 
	&\le C\lambda^{\sigma_1}
	      \|f\|_{L^{\H(p_0,p_1)}_{\theta_0+\theta_1}}
\end{align}
\end{subequations}
with an arbitrarily large $N\ge0$ and exponent $\sigma_0$,
$\sigma_0'$, $\sigma_1$ such that
\begin{equation} 
\label{sigma}
 \sigma_0>-2,
 \qquad \sigma_0'>-2,
 \qquad \sigma_1>-N+d-1.
\end{equation}
The constant $C>0$ has to be independent of $\lambda$.
Assumption \eqref{sigma} ensures that the integrals
$$\int_0^T \|F(t-s)*\left(u\otimes u\right)(s)\|_{L^{p_0}_{\theta_0}}\,ds
,\qquad
\int_0^T \|F(t-s)*\left(B\otimes B\right)(s)\|_{L^{p_0}_{\theta_0}}\,ds$$
and
$$\int_0^T \|G(t-s)*\left(u\otimes B\right)(s) \|_{L^{p_1}_{\theta_1}}\,ds$$
converge.


\section{Convolution estimates in weighted spaces}
\label{section-convolution}

The fundamental estimates \eqref{uu}-\eqref{uB}
will be a simple consequence of the following proposition.

\begin{proposition} 
\label{proposition1}
Let $a,p\in[1;+\infty]$ and $\alpha,\theta\geq0$.
For any real numbers $\lambda>0$ and $N\geq1$ let us set
$$ \Gamma_\lambda^N(x)=(\lambda+|x|)^{-N}.$$
Let also $f\in L^a_\alpha(\R^d)$ and $N>d$.
\begin{enumerate} 
\item
Then $\Gamma_\lambda^N*f\in L^p_\theta(\R^d)$, provided that
\begin{equation}\label{i)}
\theta  \le \alpha \qquad\text{and}\qquad
\theta+\frac{d}{p} \le \min\biggl\{ N-\varepsilon_{1/p}\,;
			\alpha+\frac{d}{a}-\varepsilon_{\alpha-\theta}
	 		\biggr\}.
\end{equation}
Moreover, if $N\not=d(1+\frac{1}{p}-\frac{1}{a})$, then
there exists $C>0$ such that
\begin{equation}
\label{ineg1}
\|\Gamma_\lambda^N*f\|_{L^p_\theta}\le C\lambda^{-N}(1+\lambda)^N\|f\|_{L^a_\alpha}.
\end{equation}
\item
If one assumes in addition that
\begin{equation}
\label{ii)}
\frac{1}{a}<\frac{1}{p}+\frac{1}{d},
\end{equation}
then there exists $\epsilon>0$ and two constants $C,m>0$ such that
\begin{equation}
\label{ineg2}
\|\Gamma_\lambda^N*f\|_{L^p_\theta}\le C\lambda^{-N+d-1+\epsilon}(1+\lambda)^m\|f\|_{L^a_\alpha}.
\end{equation}
\end{enumerate}
\item
When $N=d(1+\frac{1}{p}-\frac{1}{a})$, the bounds \eqref{ineg1} and \eqref{ineg2}
hold with an additional factor $(1+|\log\lambda|)$ in the right-hand sides.
In \eqref{ineg1} and \eqref{ineg2} the constant $C$
may depend on $\theta$, $a$, $\alpha$, $N$
and $d$, but it does not depend on $\lambda$ or $f$.
\end{proposition}

\begin{remark}
\begin{rm}
We shall see in the proof that we can take
\begin{gather*}
\epsilon=\min\left\{\frac{d}{p}-\frac{d}{a}+1 \semicolon \frac{N-d+1}{2}\right\},
\\
m=\max\left\{N-d+1-2\epsilon \semicolon -N+d\left(\frac{1}{p}-\frac{1}{a}+1\right)
\right\}.
\end{gather*}
\end{rm}
\end{remark}

\Proof
We start by observing that by H\"older's inequality,
\begin{equation}
\label{inj}
\|f\|_{L^q}\le C\|f\|_{L^a_\alpha} \qquad 
\text{if}\enspace \frac{1}{a}\le \frac{1}{q}\le \min\left\{1 \semicolon
\frac{1}{a}+\frac{\alpha}{d}-\varepsilon_\alpha\right\}. 
\end{equation}
Next we have
$$(1+|x|)^\theta \,|\Gamma_\lambda^N\ast f(x)|
\leq \biggl[\int_{\R^d} \Gamma_\lambda^N(x-y) |f(y)|\,dy\biggr] (1+|x|)^\theta
= I_{\theta,\lambda}(x)+J_{\theta,\lambda}(x)+K_{\theta,\lambda}(x),$$
with the following definitions :
\begin{align*}
I_{\theta,\lambda}(x) &=
\biggl(\int_{|y|\ge |x|/2} \Gamma_\lambda^N(x-y)|f(y)|\,dy\biggr)\,
(1+|x|)^\theta,\\
J_{\theta,\lambda}(x) &=
\biggl(\int_{|y|\le |x|/2} \Gamma_\lambda^N(x-y)|f(y)| \,dy\biggr)\,
(1+|x|)^\theta \:\indicatrice_{B(0,1)}(x),\\
K_{\theta,\lambda}(x) &=
 \biggl(\int_{|y|\le |x|/2 }\Gamma_\lambda^N(x-y)|f(y)| \,dy\biggr)\,
 (1+|x|)^\theta\:\indicatrice_{B(0,1)^c}(x).
\end{align*}
Here and below, $B(0,1)$ denotes the unit ball and $\indicatrice_E$ is the indicator
function of a set $E\subset\R^d$.

\paragraph{\bf The bound for $K_{\theta,\lambda}$.}
Since $|y|\le |x|/2$, we have 
$$(\lambda+|x-y|)^{-N}\le 2^N(\lambda+|x|)^{-N}.$$
Hence, using \eqref{inj} with $\frac{1}{q'}=1-\frac{1}{q}=\big(
1-\frac{\alpha}{d}-\frac{1}{a}+\varepsilon_\alpha\big)^+$,
\begin{align*}
0\le K_{\theta,\lambda}(x) 
& \leq C\,(\lambda+|x|)^{-(N-\theta)}
\int_{|y|\le \frac{|x|}{2}} |f(y)|\,dy\\
& \leq C\,(\lambda+|x|)^{-(N-\theta)} \: \|f\|_{L^q} \:
\|\indicatrice_{B(0,|x|/2)}\|_{L^{q'}}\\[2pt]
& \leq C\,(\lambda+|x|)^{-(N-\theta)} \: |x|^{[d-(\alpha+\frac{d}{a})+
	\varepsilon_\alpha]^+} \|f\|_{L^a_\alpha}.
\end{align*}
As $|x|\ge1$, it follows that $\|K_{\theta,\lambda}\|_{L^p} \leq
C\, \|f\|_{L^a_\alpha}$,
uniformly for $\lambda>0$, provided that
\begin{equation}
\label{cond K}
\theta+\frac{d}{p} \le N-\biggl[d-
\biggl(\alpha+\frac{d}{a}\biggr)+\varepsilon_\alpha\biggr]^+
-\varepsilon_{1/p}.
\end{equation}
Since $N>d$, this condition is weaker than \eqref{i)}.

\paragraph{The bound for $J_{\theta,\lambda}$.}
Using \eqref{inj} again, but with $q=a$, gives us
\begin{align*}
 0\le J_{\theta,\lambda}(x)
  &\leq C\,\indicatrice_{B(0,1)}(x)\:(\lambda+|x|)^{-N}
  \int_{|y|\le \frac{|x|}{2}} |f(y)|\,dy\\
  &\leq C\,\indicatrice_{B(0,1)}(x)\:(\lambda+|x|)^{-N}
  |x|^{d(1-1/a)} \, \|f\|_{L^a},
\end{align*}
whence
$$ \|J_{\theta,\lambda}\|_{L^p}
 \le C\,\biggl[ \lambda^{-Np}\int_{|x|\le \lambda} |x|^{dp(1-1/a)}\,dx
+ \indicatrice_{\{\lambda<1\}}
\int_{\lambda\le |x|\le 1} |x|^{-Np+dp(1-1/a)}\,dx\biggr]^{1/p}
\,\|f\|_{L^a}.$$
Thus, for all $\theta\ge0$ and $p\in[1,+\infty]$, we have
\begin{subequations}
\begin{alignat}{2}
 \|J_{\theta,\lambda}\|_{L^p}&\le
 C\left(1+\lambda^{-N+d+\frac{d}{p}-\frac{d}{a}}\right) \, \|f\|_{L^a}
 \qquad & \text{if}\enspace N\not=d\left(1+\frac{1}{p}-\frac{1}{a}\right), \\
 \text{and}\qquad
 \|J_{\theta,\lambda}\|_{L^p}&\le
 C\left(1+|\log\lambda|\right) \, \|f\|_{L^a}
 \qquad & \text{if}\enspace N=d\left(1+\frac{1}{p}-\frac{1}{a}\right).
\end{alignat}
\end{subequations}
Note that $\|J_{\theta,\lambda}\|_{L^p}$ is bounded by the right-hand side
of \eqref{ineg1}.
Moreover, if $\frac{1}{a}<\frac{1}{p}+\frac{1}{d}$, 
then $\|J_{\theta,\lambda}\|_{L^p}$ is also bounded by
the right-hand side of \eqref{ineg2},
provided that $0<\epsilon\le d(\frac{1}{p}-\frac{1}{a}+\frac{1}{d})$.

\paragraph{The bound for  $I_{\theta,\lambda}$.}
Set $F(x)=(1+|x|)^\alpha\,|f(x)|$, so that $F\in L^a(\R^d)$ and
$$0\le I_{\theta,\lambda}(x)\le
C\,(1+|x|)^{-(\alpha-\theta)}\int_{\R^d} \Gamma_\lambda^N(x-y)F(y)\,dy.$$
But
$\Gamma_\lambda^N\in L^b_\beta(\R^d)$ for all $b\in[1;+\infty]$ and
$\beta\geq0$ such that $\beta+\frac{d}{b}\le N-\varepsilon_{1/b}$.
Moreover, one has
\begin{equation}
\label{provis bound}
\|\Gamma_\lambda^N\|_{L^b_\beta}\le C
\lambda^{-N+\frac{d}{b}}(1+\lambda)^\beta.
\end{equation}

\medskip

The remaining part of the proof of Proposition~\ref{proposition1}
relies on the following lemma.

\begin{lemma}
\label{lemma1}
\begin{subequations}
Let $a,b,p\in[1;+\infty]$ and $\alpha,\beta,\theta\geq0$.
For $f\in L^a_\alpha(\R^d)$, $g\in L^b_\beta(\R^d)$, define
$$I_\theta(x)= (1+|x|)^{-(\alpha-\theta)} \, F*g (x)$$
with $F(x)=(1+|x|)^\alpha\,|f(x)|$.
If there exists $s\in[1,+\infty]$ such that:
\begin{equation}
\label{cond lemma1}
\left\{
\begin{aligned}
&\theta\leq \alpha\\
&\frac{d}{s}\leq \min\left\{
      \frac{d}{a}  \semicolon
      \left(\alpha+\frac{d}{a}\right)-\left(\theta+\frac{d}{p}\right)
        -\varepsilon_{\alpha-\theta} \semicolon
      d\left(1-\frac{1}{b}\right) 
      \right\}\\
&\frac{d}{s}\geq\max\left\{ \frac{d}{a}-\frac{d}{p} \semicolon
       \left[d-\left(\beta+\frac{d}{b}\right)+\varepsilon_\beta\right]^+
       \right\}
\end{aligned}
\right.
\end{equation}
then $I_\theta\in L^p(\R^d)$ and
\begin{equation}\label{cclemma}
\|I_\theta\|_{L^p}\le C \|f\|_{L^a_\alpha}\|g\|_{L^b_\beta}.
\end{equation}
\end{subequations}
\end{lemma}

\Proof
According to \eqref{inj}, we have $g\in L^{s'}(\R^d)$ for all $s'\in [1;+\infty]$ such that
$$\frac{1}{b}\le \frac{1}{s'}\le
\min\left\{1\semicolon\frac{1}{b}+\frac{\beta}{d}-\varepsilon_\beta\right\}.$$
Let $\frac{1}{s}+\frac{1}{s'}=1$.
We now use that $\frac{1}{a}-\frac{1}{s}\ge 0$.
The Young exponent $\Y(a,s')$ of  $a$ and $s'$ is well defined by
$\frac{1}{\Y(a,s')}=\frac{1}{a}-\frac{1}{s}$.
Moreover, one has $F*g\in L^{\Y(a,s')}(\R^d)$, {\it i.e.}
$$I_\theta \in L^{\Y(a,s')}_{\alpha-\theta}.$$
Since $\theta\le \alpha$, \eqref{inj} implies that
$I_\theta\in L^p(\R^d)$ for all $p$ such that
$$ \frac{1}{a}-\frac{1}{s} 
   \le \frac{1}{p}
   \le \min\left\{1\semicolon \frac{1}{a}-\frac{1}{s}+
     \frac{\alpha-\theta}{d}-\varepsilon_{\alpha-\theta}\right\},$$
and \eqref{cclemma} is satisfied.
\endProof

\bigskip
Let us now come back to the proof of Proposition~\ref{proposition1}.
We are going to apply the lemma with $g=\Gamma_\lambda^N$,
$I_\theta=I_{\theta,\lambda}$, $b=+\infty$ and $\beta=N$.
\begin{itemize}
\item[--]
If $\frac{1}{a}\le \frac{1}{p}$, then we further choose $s=+\infty$ and
conditions \eqref{cond lemma1} boil down (recall that $N>d$) 
to the only restriction
$\theta+\frac{d}{p}\le \alpha+\frac{d}{a}-\varepsilon_{\alpha-\theta}$.
\item[--]
If $\frac{1}{a}>\frac{1}{p}$, then we choose $\frac{1}{s}=\frac{1}{a}-\frac{1}{p}$.
In this case conditions \eqref{cond lemma1} boil down to $\theta\le \alpha$.
\end{itemize}
The first part of Proposition~\ref{proposition1} now follows
from the bounds obtained for $I_{\theta,\lambda}$, $J_{\theta,\lambda}$ and
$K_{\theta,\lambda}$.

\bigskip
To prove \eqref{ineg2}, we fix $\epsilon$ such that
$0<\epsilon\le \frac{N-d+1}{2}$.
Then we apply Lemma~\ref{lemma1} again with $g=\Gamma_\lambda^N$ and 
$I_\theta=I_{\theta,\lambda}$, but with
$b$ and  $\beta$ defined by
$$\frac{d}{b}=d-1+\epsilon, \qquad \hbox{and}\qquad \beta=N-d+1-2\epsilon.$$
By \eqref{provis bound}, one has
$\Gamma_\lambda^N\in L^b_\beta(\R^d)$ with 
$\|\Gamma_\lambda^N\|_{L^b_\beta}\le \lambda^{-N+d-1+\epsilon}\phi(\lambda)$
and $\phi\in L^\infty_{loc}([0;+\infty))$.

\bigskip\pagebreak[2]
As before,
\begin{itemize}
\item[--]
if $\frac{1}{a}\le \frac{1}{p}$, 
then we choose $s=+\infty$ in \eqref{cond lemma1}
and Lemma~\ref{lemma1} implies that 
\begin{equation} 
\label{cond I}
\|I_{\theta,\lambda}\|_{L^b_\beta}
 \le \lambda^{-N+d-1+\epsilon}\phi(\lambda)\|f\|_{L^a_\alpha},
\end{equation}
provided that
$\theta+\frac{d}{p}\le \alpha+\frac{d}{a}-\varepsilon_{\alpha-\theta}$.
\item[--]
If $\frac{1}{a}>\frac{1}{p}$, then $\frac{1}{s}=\frac{1}{a}-\frac{1}{p}$
leads again to \eqref{cond I}, provided that
$\theta\le \alpha$ and
$\frac{1}{a}\le \frac{1}{p}+\frac{1}{d}-\frac{\epsilon}{d}$.
\end{itemize}
The proof of Proposition~\ref{proposition1} is now complete.
\endProof


\section{End of the proof of Theorems~\ref{main-thm} and \ref{theorem1}}
\label{section-lifetime}

\subsection{Existence of a unique mild solution in weighted spaces}

We are now in a position to prove Theorem \ref{theorem1}.

\medskip
Under the assumptions of Theorem~\ref{theorem1}, one applies \eqref{ineg2}
with $N=d+1$ and with $\epsilon=1-\frac{d}{p_0}$ or 
$\epsilon=1-\big(\frac{2d}{p_1}-\frac{d}{p_0}\big)^+$ respectively;
assumption \eqref{ii)} is ensured by \eqref{hypo}.
This proves \eqref{uu} and \eqref{BB} with
$$\sigma_0=-1-\frac{d}{p_0} \qquad\text{and}\qquad
\sigma_0'=-1-\left(\frac{2d}{p_1}-\frac{d}{p_0}\right)^+.$$

A new application of \eqref{ineg2}
with any $N$ such that $N\ge \max\{d+1 \semicolon
\theta_1+\frac{d}{p_1}\}+\varepsilon_{1/p_1}$ 
and $\epsilon=1-\frac{d}{p_0}$ yields \eqref{uB} with
$\sigma_1=-N+d-d/p_0$.

With the preceeding values of $\sigma_0,\sigma_0'$ and $\sigma_1$,
the assumption \eqref{hypo} implies \eqref{sigma}.
As indicated in section~\ref{section-generalities}, this yields \eqref{bicontinuity}
and ensures that the operator norm of~$\V$ tends to zero as a power of $T$,
when $T\to0$ :
$$\interleave\V\interleave_{{\cal C}([0,T];X)} \leq C \,
\max\left\{ T^{1+\frac{\sigma_0}{2}} \,;\,  T^{1+\frac{\sigma_0'}{2}} \,;\,
T^{1+\frac{1}{2}(\sigma_1+N-d-1)} \right\}.$$
This ensures finally the conclusions
\eqref{mild-solution_Lp} and \eqref{mild-solution}
of Theorem~\ref{theorem1}.

More precisely, our argument proves that under the assumptions
of Theorem~\ref{theorem1}, the maximal lifetime~$T^*$ of the
mild solution in $X=L^{p_0}_{\theta_0}\times L^{p_1}_{\theta_1}$ satisfies
\begin{equation}
\label{To deduce...}
 T^\ast \geq c\,\min\left\{ 1 \semicolon
   \|(u_0,B_0)\|_X^{ -2/(1-\frac{d}{p_0}) }\, \semicolon
   \|(u_0,B_0)\|_X^{ -2/\big(1-\big[\frac{2d}{p_1}-\frac{d}{p_0}\big]^+
      \big) }
   \right\},
\end{equation}
with a constant $c>0$, depending on all the parameters, but
not on $u_0$ or on $B_0$.

\subsection{Comparison of lifetimes in Theorem~\ref{theorem1}}
It only remains to establish that lifetimes are independent of the
admissible pairs of indices chosen to construct the solution.

\begin{proposition}
\label{proposition2}
Let $u_0\in L^{p_0}_{\theta_0}(\R^d)\cap L^{\tilde{p_0}}_{\tilde\theta_0}(\R^d)$ and
$B_0\in L^{p_1}_{\theta_1}(\R^d)$.
Set $\eta_0=\theta_0+d/p_0$, $\tilde{\eta_0}=\tilde\theta_0+d/\tilde{p_0}$
and $\eta_1=\theta_1+d/p_1$.
Assume that $d\geq2$ and
\begin{equation}
\label{parameters}
\left\{\begin{aligned}
&d<p_0,\tilde{p_0}\le +\infty\\
&\frac{2}{p_1}< \min\left\{ \frac{1}{p_0}+\frac{1}{d} \semicolon
          \frac{1}{\tilde{p_0}}+\frac{1}{d} \right\}\\
&\eta_0\le \min\left\{d+1-\varepsilon_{1/p_0} \semicolon
          2\eta_1-\varepsilon_{2\theta_1-\theta_0} \semicolon
          2\eta_1+\frac{d}{p_0}-\frac{2d}{p_1}
	  \right\} \\
&\tilde{\eta_0} \le \min\left\{d+1-\varepsilon_{1/\tilde{p_0}} \semicolon
          2\eta_1-\varepsilon_{2\theta_1-\tilde\theta_0} \semicolon
	  2\eta_1+\frac{d}{\tilde{p_0}}-\frac{2d}{p_1}
	  \right\}.
\end{aligned}\right.
\end{equation}
Let $T^\ast$ and $\tilde{T}$ be the lifetimes of the solution $(u,B)$ of (mhd) emanating from $(u_0,B_0)$
in the respective weighted spaces, {\it i.e.}
\begin{align*}
T^\ast =&\sup\left\{\, T>0 \enspace\text{s.t.}\enspace
  (u,B)\in {\cal C}([0,T];L^{p_0}_{\theta_0}\times L^{p_1}_{\theta_1})
  \,\right\},\\
\tilde{T} =&\sup\left\{\, T>0 \enspace\text{s.t.}\enspace
  (u,B)\in {\cal C}([0,T];L^{\tilde{p_0}}_{\tilde\theta_0}\times
  L^{p_1}_{\theta_1}) \, \right\}.
\end{align*}
Then $\tilde{T}=T^\ast$.
\end{proposition}

\Proof
The structure of the proof is similar to that of \cite{Vig05}.
Let us assume that we have, for example,~$\tilde T<T^*$.
Unicity of mild solutions ensures that they agree on
$[0,\tilde T[$. We are going to prove that
$$\sup_{t\in[0,\tilde T[} \left( \|u(t)\|_{L^{\tilde p_0}_{\tilde \theta_0}}
+ \|B(t)\|_{L^{p_1}_{\theta_1}} \right)<+\infty.$$
Then \eqref{To deduce...} would imply that the mild solution~$(u,B)$
in~${L^{\tilde{p_0}}_{\tilde \theta_0}\times L^{p_1}_{\theta_1}}$
could be extended beyond $\tilde T$, and that would contradict the definition of
$\tilde T$.

\medskip
First of all, let us recall (see, {\it e.g.}, \cite[\S2.2]{Vig05}) that
there exists a constant $C_0>0$ depending only on $d$ and $\theta$,
such that
\begin{equation}    
\sup_{\tau\in [0,\tilde T]}\|e^{\tau\Delta}v\|_{L^{p_1}_{\theta_1}}
   \le C_0\,(1+\tilde T)^{\theta_1/2} \, \|v\|_{L^{p_1}_{\theta_1}}.
\end{equation}
In the following, we set $A = C_0\,(1+\tilde T)^{\theta_1/2}$.

Note also that we can obviously assume that $u\not\equiv0$ in $[0,\tilde T]$.

\paragraph{The bound for $B$.}
By the second of the integral equations (IE), one has for $0\le s\le t< \tilde T$~:
$$B(t)=e^{(t-s)\Delta}B(s)-\int_s^t G(t-\tau)*\left(u\otimes B-B\otimes u\right)(\tau)\,d\tau.$$
Proposition~\ref{proposition1} applied to the upper
bound of $G$ given by \eqref{pointwise}, with $\epsilon=1-\frac{d}{p_0}$
in \eqref{ineg2}, yields
$$\forall\;\tau\le t\le \tilde T,\qquad
\|G(t-\tau)*\left (u\otimes B\right)(\tau)\|_{L^{p_1}_{\theta_1}}
  \le K(t-\tau)^{-\sigma} \| \left(u\otimes B\right)(\tau)\|_{L^{\H(p_0,p_1)}_{\theta_0+\theta_1}}$$
where $\sigma=\frac{1}{2}(1+\frac{d}{p_0})$ and $K$ is a constant, possibly
depending on $T^\ast$ and all the parameters contained in \eqref{parameters},
but not on $\tilde T$.
Note that $\sigma<1$.
Thus, for all $t\in[0;\tilde T]$,
\begin{equation}
\label{norme B}
\|B(t)\|_{L^{p_1}_{\theta_1}}\le 
	A\,\|B(s)\|_{L^{p_1}_{\theta_1}} +
	 K\,\frac{ (t-s)^{1-\sigma}}{1-\sigma } \:
            \sup_{\tau\in[s,t]} \|u(\tau)\|_{L^{p_0}_{\theta_0}} \cdot
	    \sup_{\tau\in[s,t]} \|B(\tau)\|_{L^{p_1}_{\theta_1}}  \,.
\end{equation}
Now let $(T_n)_{n\ge0}$ be the increasing sequence defined by
$$T_n = n \Delta \qquad\text{with}\qquad
\Delta = \left(\frac{2 K}{1-\sigma}\,
  \sup_{\tau\in [0,\tilde T]}\|u(\tau)\|_{L^{p_0}_{\theta_0}}
\right)^{-1/(1-\sigma)}$$
and $N\in\mathbb{N}$ such that $T_N\le \tilde T< T_{N+1}$.
For $0\le n\le N$, let $I_n$ be the interval $[T_n,T_{n+1}]\cap [0,\tilde T[$
and  
$$ M_n = \sup_{\tau\in I_n}\|B(\tau)\|_{L^{p_1}_{\theta_1}}.$$
Applying \eqref{norme B} with $s=T_n$ and $t\in I_n$ for $n=0,\dots, N$, we get
$$ M_0\le 2A\|B_0\|_{L^{p_1}_{\theta_1}}
   \qquad\text{and}\qquad
   M_n\le 2AM_{n-1}
   \quad(1\le n\le N),$$
whence
$$\sup_{t\in[0,\tilde T[} \|B(t)\|_{L^{p_1}_{\theta_1}}=\max_{0\le n\le N} M_n
\leq (2A)^{N+1}\|B_0\|_{L^{p_1}_{\theta_1}}.$$
Finally, this leads to~:
\begin{equation}\label{B borne}
\sup_{t\in[0,\tilde T[} \|B(t)\|_{L^{p_1}_{\theta_1}}
\leq C\,\|B_0\|_{L^{p_1}_{\theta_1}}
   \exp\left(
   \biggl(1+\tilde T \sup_{s\in[0,\tilde T]}\,
      \|u(s)\|_{L^{p_0}_{\theta_0}}^{2/(1-\frac{d}{p_0})} \biggr)
   \left(1+\theta_1\log(1+\tilde T)\right)
   \right).
\end{equation}
The right-hand side is finite because we assumed $\tilde T<T^\ast$.

\paragraph{The bound for $u$.}
For $0\le s\le t< \tilde T$, one has
$$ u(t)=e^{(t-s)\Delta}u(s)
  -\int_s^t F(t-\tau)\ast (u\otimes u)(\tau)\,d\tau
  +\int_s^t F(t-\tau)\ast (B\otimes B)(\tau)\,d\tau.$$
Proposition~\ref{proposition1}, applied this time to the upper
bound of $F$ given by \eqref{pointwise}, yields
\begin{align*}
\|u(t)\|_{L^{\tilde p_0}_{\tilde\theta_0}} &\le
	A\|u(s)\|_{L^{\tilde p_0}_{\tilde\theta_0}} +
	K\, \frac{(t-s)^{1-\sigma}}{1-\sigma}
		\sup_{\tau\in[s,t]} \|u(\tau)\|_{L^{p_0}_{\theta_0}} 
		\cdot
		\sup_{\tau\in[s,t]} \|u(\tau)\|_{L^{\tilde p_0}_{\tilde \theta_0}}\\
      &  \qquad\qquad
       + K \, \frac{ (t-s)^{1-\tilde \sigma}}{1-\tilde \sigma }
		\left( \sup_{\tau\in[s,t]} \|B(\tau)\|_{L^{p_1}_{\theta_1}}
		\right)^2
\end{align*}
with $\sigma=\frac{1}{2}(1+\frac{d}{p_0})$ and
$\tilde\sigma=\frac{1}{2}(1+(\frac{2d}{p_1}-\frac{d}{p_0})^+)$.
Note that $\sigma$ is the same as before and that $\tilde \sigma<1$;
$K$ depends on $T^\ast$ and all the parameters, except $\tilde T$.
The last term is uniformly bounded by
$$L=\frac{ K \, \tilde T^{1-\tilde \sigma}}{1-\tilde \sigma}
    \left( \sup_{\tau\in[0,\tilde T[} \|B(\tau)\|_{L^{p_1}_{\theta_1}} \right)^2$$
which is a finite constant because \eqref{B borne} holds.
Define $(T_n)_{n\ge0}$ and $I_n$ as before. Let also
$$\tilde M_n=\sup_{\tau\in I_n} \|u(\tau)\|_{L^{\tilde p_0}_{\tilde
\theta_0}}.$$
Recall that $N$ is the integer part of $\tilde T/\Delta$.
Then, for $1\le i\le N$, one has
$$\tilde M_0\le 2A\,\|u_0\|_{L^{\tilde p_0}_{\tilde\theta_0}}+2L
  \qquad\text{and}\qquad
  \tilde M_n\le 2A\,\tilde M_{n-1}+2L,$$
hence
$$\sup_{t\in[0,\tilde T[} \|u(t)\|_{L^{\tilde{p_0}}_{\tilde \theta_0}}
=\max_{0\le n\le N} \tilde M_n
\le (2A)^{N+1}\,\|u_0\|_{L^{\tilde p_0}_{\tilde \theta_0} }
+2L \left[1+\ldots+(2A)^{N-1}+(2A)^N\right] < +\infty.$$
Combined with \eqref{To deduce...} and~\eqref{B borne},
this estimate ensures that $\tilde T \geq T^\ast$.
Exchanging the roles of~$\tilde T$ and~$T^\ast$, one finally obtains
that $\tilde T = T^\ast$.
\endProof

\medskip
An analogous result holds if we assume instead 
$u_0\in L^{p_0}_{\theta_0}(\R^d)$ and  
$B_0\in L^{p_1}_{\theta_1}(\R^d)\cap L^{\tilde p_1}_{\tilde \theta_1}(\R^d)$,
with obvious modifications in \eqref{parameters}~:
\begin{equation}\tag{\ref{parameters}'}
\left\{\begin{aligned}
&d<p_0\le +\infty\\
&\max\left\{\frac{2}{p_1} \semicolon \frac{2}{\tilde{p_1}} \right\}
          < \frac{1}{p_0}+\frac{1}{d} \\
&\eta_0 \le 
          \min\left\{d+1-\varepsilon_{1/p_0} \semicolon
          2\eta_1-\varepsilon_{2\theta_1-\theta_0} \semicolon
	  2\tilde{\eta_1}-\varepsilon_{2\theta_1-\theta_0}
	  \right\}.\\
&\eta_0 \le \min\left\{
          2\eta_1+\frac{d}{p_0}-\frac{2d}{p_1} \semicolon
          2\tilde{\eta_1}+\frac{d}{p_0}-\frac{2d}{\tilde{p_1}}
	  \right\}.
\end{aligned}\right.
\end{equation}
Theorem~\ref{theorem1} is now established.

\subsection{The proof of Theorem~\ref{main-thm}}
\label{proof-main-thm}
Let $p_0,p_1$ and $\theta_0,\theta_1$ such that \eqref{indices-thm} and
\eqref{eta-thm} hold.

\medskip
If $\theta_0\le 2\theta_1$, $p_0\leq d/\delta-\varepsilon_\delta$
and $\eta_0 \leq d+1-\varepsilon_{1/p_0}$, then \eqref{hypo} and \eqref{local}
hold, and there is nothing more to prove since Theorem~\ref{theorem1} already
gives a stronger conclusion.

\medskip
In all the other cases and for any $\epsilon>0$, our assumptions yield
an embedding $L^{p_0}_{\theta_0}\subset L^{q}_{\mu}$
such that Theorem~\ref{theorem1} may be applied to
$$(u_0,B_0)\in L^q_\mu\times L^{p_1}_{\theta_1}$$
and with
$$\mu+\frac{d}{q}=\eta_0-\epsilon.$$
It follows that $u\stackrel{L^2}{=}{\cal O}(|x|^{-(\eta_0-\epsilon)})$
and $B\stackrel{L^2}{=}{\cal O}(|x|^{-\eta_1})$  when $|x|\to +\infty$.
Letting $\epsilon\to0$, this will conclude the proof of Theorem~\ref{main-thm}.

\medskip
Let us be more precise about the embedding
$L^{p_0}_{\theta_0}\subset L^{q}_{\mu}$.
Actually, various choices are possible for $(q,\mu)$. We have chosen
the indices that are represented on the interpolation diagram
(see Fig.~1 p.~\pageref{figure}) by a {\sl dash-dotted} line. 

If the magnetic field decays sufficiently fast, namely if $\eta_1\geq
(d+1+\delta)/2$, the only case not included in Theorem~\ref{theorem1} is
that of $\eta_0=d+1$ with $p_0$ finite. In this case, one may take
$$(q,\mu)=(p_0,\theta_0-\epsilon).$$

Let us now assume that $\eta_1\leq (d+1+\delta)/2$ and, for the moment,
that $p_1\geq 2d$. Then the cases to be dealt with correspond either to
$\theta_0>2\theta_1$ or to $\eta_0=2\eta_1$, or to both.

-- If $\theta_0>2\theta_1$, then
$$\frac{d}{q}=\theta_0 -2\theta_1+\frac{d}{p_0}-\epsilon
\qquad\text{and}\qquad\mu=2\theta_1$$
are suitable, even if $\eta_0=2\eta_1$.

-- If $\theta_0\leq2\theta_1$ and $\eta_0=2\eta_1$, one may again choose
$(q,\mu)=(p_0,\theta_0-\epsilon)$.
 
\smallskip
Finally, if $d<p_1<2d$ and $\eta_1\leq (d+1+\delta)/2$, one may use the
following barrier :
$$\frac{d}{q} = 1-(1-\delta)\kappa
\qquad
\mu=2\theta_1(1-\kappa)
\quad\text{and}\quad
\kappa=1-\frac{\eta_0-\delta-\epsilon}{2(\eta_1-\delta)}\cdotp$$
The proof of Theorem~\ref{main-thm} is now complete.
\endProof


\section{Instantaneous spreading of rapidly decreasing fields}
\label{section-spreading}

This section is included for completeness and contains the proof
of theorem \ref{thm_optimal}, and some remarks about exceptional
solutions to (MHD) that decay extremely fast.

\subsection{Proof of theorem \ref{thm_optimal}}

Following \cite{BraM02}, we define $E$ as the space of all functions $f\in
L^1_{\textrm{loc}}(\R^d)$ such that
\begin{equation}
\label{norm-E}
\|f\|_E \definition \int_{|x|\le1}|f(x)|\,dx
\,+\, \sup_{R\ge1} \: R\int_{|x|\ge R}|f(x)|\,dx
\end{equation}
is finite, and
$$\lim_{R\to +\infty} R\int_{|x|\ge R} |f(x)|\,dx=0.$$
H\"{o}lder inequality implies that :
$$L^{p_0}_{\theta_0}(\R^d)\subset E \qquad\text{whenever}\qquad
 \left\{\begin{array}{l}\displaystyle
    \theta_0+\frac{d}{p_0}\ge d+1 \quad (p_0<+\infty) \quad
    \text{or}\\[1em]
    \theta_0>d+1 \quad (p_0=+\infty).
\end{array}
\right.$$
Let us prove that $\|u\|_E$ cannot remain uniformly bounded during a
positive time interval, unless the orthogonality
relations \eqref{orthogonality} are satisfied.

\begin{proposition}
\label{theorem2}
Let $(u,B)\in {\cal C}([0,T];L^2(\R^d)\times L^2(\R^d))$ a solution to (MHD)
such that $u_0\in E$. Assume that
\begin{subequations}
\label{localization}
\begin{gather}
\label{localization1}
u\in L^\infty([0,T];E)\\
\label{localization2}
|u|^2 + |B|^2 \in L^\infty([0,T];E).
\end{gather}
\end{subequations}
Then there exists a constant $c\ge0$ such that the components of
the initial data satisfy
\begin{equation}
\label{orth}
\forall j,k \in\{1,\ldots,d\}, \qquad
\int_{\R^d} u_0^j u_0^k-B_0^j B_0^k=c\,\delta_{j,k},
\end{equation}
where $\delta_{j,k}=1$ if $j=k$ and $0$ otherwise.
\end{proposition}

\Proof
The proof will only be sketched briefly since it is a straightforward
adaptation  of \cite{BraM02}.
Let us write the first equation of (MHD) in the following form
(recall that $S$ and $R_e$ can be set equal to~1):
\begin{equation} 
\label{inteq}
u(t)-e^{t\Delta}u_0
+\sum_{j=1}^d\int_0^t e^{(t-s)\Delta}\,\partial_j(u^ju-B^jB)\,ds=
-\int_0^t e^{(t-s)\Delta} \, \nabla P(s)\,ds,
\end{equation}
where $P=p+\frac{|B|^2}{2}$ is the total pressure.
Arguing as in \cite{BraM02}, we see that \eqref{localization} imply
that all the terms in the left-hand side of \eqref{inteq}
belong to $L^\infty([0,T];E)$.
Thus, we have
$$\nabla \tilde P\in L^\infty([0,T];E)
\qquad\text{with}\qquad
\tilde P(t)=\int_0^t e^{(t-s)\Delta}\,P(s)\,ds.$$
Let 
$$\tilde u^{j,k}(t)=\int_0^t e^{(t-s)\Delta} \, u^ju^k(s)\,ds$$
and
$$\tilde B^{j,k}(t)=\int_0^t e^{(t-s)\Delta} \, B^jB^k(s)\,ds.$$
Taking the divergence in~\eqref{inteq} yields
$$-\Delta\tilde P=\sum_{j,k=1}^d\partial_j\partial_k(\tilde u^{j,k}-\tilde
B^{j,k}).$$
One now deduces \eqref{orth},
applying Lemma~2.3 and Proposition~2.4 of \cite{BraM02}.
\endProof

\bigskip

The proof of Theorem~\ref{thm_optimal} is now very easy.
Thanks to \eqref{inclusion} and \eqref{def2_eta}, assumptions \eqref{optimality_u} and
\eqref{optimality_B} imply the existence of $\varepsilon'>\varepsilon''>0$
such that
$$\sup_{t\in[0,T]} |u(t,\cdot)|\in L^2_{\frac{d}{2}+1+\varepsilon'} \subset
L^{1}_{1+\varepsilon''}\subset E.$$
Moreover, the definition of the \emph{$L^2$ decay rate at
infinity} \eqref{distr-decay} implies that
$$\lim_{R\to\infty} R^{d+2+2\varepsilon'}
\int_{R\leq |x|\leq 2R} |u(t,x)|^2 \,dx = 0$$
and $$\lim_{R\to\infty} R^{1+\varepsilon'}
\int_{R\leq |x|\leq 2R} |B(t,x)|^2 \,dx = 0,$$
uniformly for $t\in[0,T]$. Therefore
$$\sup_{t\in[0,T]} \left(|u(t,\cdot)|^2+|B(t,\cdot)|^2\right)
\in L^1_{1+\varepsilon''} \subset E.$$
Conclusion \eqref{orthogonality} now follows from proposition~\ref{theorem2}.

\subsection{Solutions of (MHD) with an exceptional spatial behavior}

We finally observe that solutions that decay faster than predicted 
by Theorem~\ref{theorem1} do exist.

\medskip
Such solutions can be constructed starting with properly symmetric
initial data.
Assume, {\it e.g.\/}, that $u_0$ and $B_0$ are rapidly decreasing in the usual
pointwise sense when $|x|\to +\infty$
(faster than any inverse polynomial)
and that
$Au_0(x)=u_0(Ax)$, $AB_0(x)=B_0(Ax)$ for all $x\in\R^d$ and all matrix $A\in G$,
where $G$ is a subgroup of the orthogonal group $O(d)$.
Then the solution of (MHD) will inherit this property as far as it exists,
the system being invariant under rotations.
If the group $G$ is rich enough, then these symmetry relations ensure the validity
of conditions \eqref{orthogonality}.
Moreover the  decay rate of the velocity field of the corresponding solution  
will depend on the symmetry group to which $(u_0,B_0)$ belongs.

In dimension $d=2,3$ and for the Navier--Stokes equations, the optimal decay
rates of the solution have been computed in \cite{Bra04iii} for each symmetry
group.
With simple modifications in the proofs, one could show that
{\it the same} decay rates hold for the solution of (MHD).
This is not surprising: indeed, since the magnetic field decays fast when
$|x|\to +\infty$, the decay of the velocity field is governed only by
the decay rate of the kernels $F_{j,h}^k$, defined by \eqref{symbol},
and by the possible corresponding cancellations.
These kernels are the same ones that appear in the Navier--Stokes system
as well.

Thus, for example, in dimension $d=2$ and when $G$ is the cyclic group of order $n$,
one has
$$\forall t\in [0,T^\ast), \qquad
u(t,x)={\cal O}(|x|^{-(n+1)})$$ in the usual pointwise sense, when $|x|\to +\infty$.
In particular, the property of being simultaneously completely invariant under 
rotations ({\it i.e.} $G=SO(2)$) and rapidly decreasing at infinity
will be conserved by~$(u,B)$ during the evolution, if such property
already holds for $(u_0,B_0)$.

In dimension three, the largest decay rates of the velocity field
({\it i.e.} like $|x|^{-8}$ as $x\to +\infty$)
are obtained with the symmetry groups of the icosahedron.
Those symmetric solutions are however unstable: in general,
the velocity field of an infinitesimal perturbation of a highly
symmetric flow will decay much more slowly at infinity.

\bibliographystyle{amsplain}

\bigskip
\makeadress

\end{document}